\documentclass[12pt]{article}
 \usepackage{latexsym,amsfonts,amssymb}
 \usepackage{amsmath}
 \usepackage[dvips]{color}
 \usepackage{colordvi,multicol}
 \usepackage[marginal]{footmisc}

 \def \cal{\mathcal}
 
 \newtheorem{thm}{Theorem}[section]

 \newtheorem{lem}[thm]{Lemma}
 
 \newtheorem{defi}[thm]{Definition}
 
 \newtheorem{rem}[thm]{Remark}
 \newtheorem{exa}[thm]{Example}

 \numberwithin{equation}{section}
 \date{}
 
\begin{document}
 	
 	\title{\bf Random periodic solutions of non-autonomous stochastic differential equations}
 	\author{}
 	\maketitle
 	\centerline{Zhao Dong$^{1,2}$, Weili Zhang$^{1,2,*}$ and Zuohuan Zheng $^{3,1,2}$} \centerline{\small $^1$  Academy of Mathematics and Systems Science, Chinese Academy of Sciences, Beijing, 100190, China} \centerline{\small $^2$ School of Mathematical Sciences, University of Chinese Academy of Sciences, Beijing, 100049, China} \centerline{\small $^3$ College of Mathematics and Statistics, Hainan Normal University, Haikou, Hainan 571158, China}
 	\centerline{\small dzhao@amt.ac.cn, zhangweili@amss.ac.cn, zhzheng@amt.ac.cn}
   \footnote{* Corresponding author}
 	\vskip 1cm  
 	\vskip 0.5cm \noindent{\bf Abstract:}\quad In this paper, we study the existence of random periodic solutions for nonlinear stochastic differential equations with additive white noise. We extend the input-to-state characteristic operator of the system to the non-autonomous stochastic differential equation via the pull-back of the discretised stochastic differential equation. We then use the completeness of the measurable function space which we construct skillfully and the Banach fixed point theorem to prove the existence of a fixed point of the gain operator. And we prove that the image for the input-to-state characteristic operator at this fixed point is a random periodic solution for the forward stochastic flow generated by the non-autonomous stochastic differential equation. Finally, we present some examples.
 	\smallskip
 	
 	\vskip 0.5cm
 	
 	
 	\vskip 0.3cm
 	
 	
 	\section{Introduction}	
 	\quad Periodic solution has been a central concept of the qualitative theory of ordinary differential equations and the deterministic dynamic system theory, it originated from Poincar\'e's pioneering work ({\cite{P}}). The deterministic autonomous dynamical system does not depend on the initial time, but depends on the time interval. However, the remarkable feature of deterministic non-autonomous dynamical system is that its evolution depends on both the current time and the initial time. Therefore, this kind of systems require a two-parameter semi-group to characterize, which makes it describe more natural phenomena. Denote $\triangle:=\{(t,s)\in\mathbb{R}^2,s\leq t\}$. For a deterministic non-autonomous dynamical system $u(t,s): X\rightarrow X$ over time $(t,s)\in\triangle$, where $X$ is a metric space. A periodic solution for $u$ is a periodic function $y:\mathbb{R}\rightarrow X$ with periodic $T>0$ such that 
 	$$u(t,s)y(s)=y(t)\ \ \ \text{and} \ \ \ y(s+T)=y(s) \ \ \ \text{for}  \ \text{all}  \ (t,s)\in\triangle.$$
 	
 	As the random counterpart of periodic solution, H.Z. Zhao and Z-H. Zheng first proposed the definition of random periodic solution for a $C^1$-cocycle ({\cite{ZZ}}). Later the concept of random periodic solutions for semi-flows generated by non-autonomous stochastic differential equations and stochastic partial differential equations were given in {\cite{FZZ}} and {\cite{FZ}}, respectively. Also, they led to more progress on investigations of various issues in autonomous stochastic differential equations and non-autonomous stochastic differential equations. They include numerical analysis of random periodic solutions and periodic measures of stochastic differential equations ({\cite{FLZ}}, {\cite{FLZZ}}); anticipating random periodic solutions of stochastic differential equations ({\cite{FWZ}}); random periodic processes, periodic measures and ergodicity ({\cite{FFZ}}) etc. 
 	
 	
 Let $(\Omega,{\cal F},P)$ be a probability space, $X$ be a Polish space and ${\cal B}(X)$ be its Borel $\sigma$-algebra.  Denote by $(\Omega, {\cal F}, P, (\theta_s)_{s\in\mathbb{R}})$ a metric dynamical system and $\theta_s:\Omega\rightarrow\Omega$ is assumed to be measurably invertible for all $s\in\mathbb{R}$.  	
 	
 	\begin{defi}[{\cite{ZZ}}]
 		A random periodic solution of periodic $T>0$ for the random dynamical system $\varphi:\mathbb{R}\times\Omega\times X\rightarrow X$ is an ${\cal F}$-measurable map $Y: \mathbb{R}\times\Omega\rightarrow X$ such that for almost all $\omega\in\Omega$,
 		\begin {equation}\label{11}
 		\varphi(t,\theta_s\omega)Y(s,\omega)=Y(t+s,\omega),\ \ Y(s+T,\omega)=Y(s,\theta_T\omega), \ \ \rm{for} \ \rm{any} \ \ t, s\in\mathbb{R}.
 	\end{equation}	
 \end{defi}	

 Consider a forward stochastic flow  $\varphi:\triangle\times\Omega\times X\rightarrow X$ of periodic $T$, which satisfies the following standard condition
 \begin {equation}\label{12}
 \varphi(s,s)=\rm{Id},\ \ \varphi(t,s,\omega)=\varphi(t,r,\omega)\circ\varphi(r,s,\omega), 
\end{equation}
and the periodic property
\begin {equation}\label{13}
\varphi(t+T,s+T,\omega)=\varphi(t,s,\theta_T\omega),
\end{equation}
for all $s\leq r\leq t$, $s,r,t\in\mathbb{R}$. 

\begin{defi}[{\cite{FZZ}}]
	A random periodic solution of periodic $T>0$ for the forward stochastic flow $\varphi:\triangle\times\Omega\times X\rightarrow X$ is an ${\cal F}$-measurable map $Y: \mathbb{R}\times\Omega\rightarrow X$ such that for almost all $\omega\in\Omega$,
	\begin {equation}\label{14}
	\varphi(t,s,\omega)Y(s,\omega)=Y(t,\omega), Y(s+T,\omega)=Y(s,\theta_T\omega), \ \ \rm{for} \ \rm{any} \ \ (t, s)\in\triangle.
\end{equation}	
\end{defi}

In this paper, we consider the following nonlinear stochastic system:
\begin{equation}\label{15}
dX_t=[AX_t+h(t,X_t)]dt+\sigma(t) dW_t, 
\end{equation}       	
where $W_t(\omega)=(W_t^1(\omega),...,W_t^m(\omega))$ is a two-side time Wiener process with values in $\mathbb{R}^m$ on the canonical Wiener space $(\Omega,{\cal F},({\cal F}^t)_{t\in\mathbb{R}},P)$, i.e., ${\cal F}$ is the Borel $\sigma$-algebra of $\Omega=C_0(\mathbb{R},\mathbb{R}^m)=\{\omega:\omega(t)$ continuous, $\omega(0)=0,t\in\mathbb{R}\}$; ${\cal F}_s^t$ is the least complete $\sigma$-field for which all $W_u-W_v$, $s\leq v\leq u\leq t$ are measurable and ${\cal F}^t={\cal F}_{-\infty}^t=\bigvee_{s\leq t}{\cal F}_s^t$; $P$ is the Wiener measure. $A=(a_{ij})_{d\times d}$ is a $(d\times d)$-dimensional matrix. $h:\mathbb{R}\times \mathbb{R}^d\rightarrow \mathbb{R}_+^d$, $h(t+T,X_t)=h(t,X_t)$, for any $t\in\mathbb{R}, x\in\mathbb{R}^d$, $T>0$ is a constant. $\sigma(t)=\sigma(t+T)$ is a $(d\times m)$-dimensional matrix.

In {\cite{JL}}, J.F. Jiang and X. Lv studied small-gain results for nonlinear stochastic system driven by additive white noise of the form with $h(t,x)\equiv h(x), \sigma(t)\equiv \sigma$ in (\ref{15}), and $h$ is order-preserving or anti-order-preserving. Motivated by their work, we extend their results to non-autonomous stochastic differential equation (\ref{15}). In our case, the solution of (\ref{15}) is a forward stochastic flow which leads to the input-to-state characteristic operator ${\cal K}(u)$ not only depending on $\omega\in\Omega$ but also depending on $t\in \mathbb{R}$. Thus, in order to get the fixed point of the gain operator, we construct a measurable function space skillfully. And we prove that the image for the input-to-state characteristic operator at this fixed point is a random periodic solution for the non-autonomous stochastic differential equation (\ref{15}). 


The rest of this paper is organized as follows. In section 2, we review some preliminary concepts and definitions, present the assumptions for the non-autonomous stochastic differential equation (\ref{15}), and define an operator of the system via the pull-back of the discretised stochastic differential equation. In section 3, we  establish some auxiliary lemmas and proposition, and present the definition of another operator and its properties. In section 4, the main result Theorem 4.2 is proved. In section 5, we present some examples.

\section{Preliminaries and Assumptions}   
\quad Before starting our main results, we introduce some basic concepts and notations.     	

\begin{defi}[{\cite{A}}]
	A family of mapping on the sample space $\Omega$, $\theta_t:\Omega\rightarrow\Omega, t\in \mathbb{R}$, is called a measurable dynamical system if the following conditions are satisfied \\	
	\quad (i) Identity property: $\theta_0=\rm{Id}$;\\	
	\quad (ii) Flow property: $\theta_{t+s}=\theta_{t}\circ \theta_{s}$;\\
	\quad (iii) Measurability: $(\omega, t)\rightarrow \theta_{t}\omega$ is measurable.\\	
	It is called a measure-preserving or metric dynamical system if, furthermore\\
	\quad (iv) Measure-preserving preserving property: $P(\theta_{t}(A))=P(A)$, for every $A\in {\cal F}$ and $t\in \mathbb{R}$.\\	
	In this case, $P$ is called an invariant measure with respect to the dynamical system $\theta_t$.
\end{defi}

\begin{defi}[{\cite{A}}]
	A measurable random dynamical system (RDS) on the measurable space $(X,{\cal B}(X))$ over a metric dynamical system $(\Omega,{\cal F},P,(\theta_{t})_{t\in\mathbb{R}})$ with time $\mathbb{R}$ is a mapping $$\Phi:{\mathbb{R}}\times\Omega\times X \rightarrow X, \  (t,\omega,x)\mapsto\Phi(t,\omega,x)$$ 
	with the following properties:\\
	\quad (i) Measurability: $\Phi$ is $({\cal B}({\mathbb R})\otimes {\cal F}\otimes {\cal B}(X), {\cal B}(X))$-measurable.\\	
	\quad (ii) Cocycle property: The mappings $\Phi(t,\omega):=\Phi(t,\omega,\cdot) : X\rightarrow X$ form a cocycle over $\theta(\cdot)$, i.e. they satisfy
	$$
	\Phi(0,\omega)={id_X},  \ \ \text{for} \ \text{all} \  \omega\in \Omega,
	$$
	$$
	\Phi(t+s,\omega)=\Phi(t,\theta_s\omega)\circ\Phi(s,\omega),  \ \ \text{for} \ \text{all} \  s,t\in\mathbb{R}, \ \omega\in \Omega.
	$$	
	Here $\circ$ means composition of mappings.
\end{defi}

\begin{defi}[{\cite{KH}}]
	Let $\varphi(t,s,x,\omega), s,t\in\mathbb{R}, x\in\mathbb{R}^d$ be a continuous $\mathbb{R}^d$-valued random field defined on the probability space $(\Omega,{\cal F},P)$. Then for almost all $\omega$, $\varphi(t,s,\omega)\equiv \varphi(t,s,\cdot,\omega)$ defines a continuous map from $\mathbb{R}^d$ into itself for any $s,t$. It is called a stochastic flow if there exists a null set $N$ of $\Omega$ such that for any $\omega\in N^c$, the family of continuous maps $\{\varphi(t,s,\omega):s,t\in\mathbb{R}\}$ defines a flow i.e. it satisfies the following properties:\\
	\quad (i) $\varphi(u,s,\omega)=\varphi(u,t,\omega)\circ\varphi(t,s,\omega)$ holds for all $s,t,u$, where $\circ$ denotes the composition of maps.\\	
	\quad (ii) $\varphi(s,s,\omega)=$ identity map for all $s$.\\
	Further if $\varphi(t,s,\omega)$ satisfies (iii), it is called a stochastic flow of homeomorphisms.\\
	\quad (iii) the map $\varphi(t,s,\omega):\mathbb{R}^d\rightarrow\mathbb{R}^d$ is an onto homeomorphism for all $s,t\in\mathbb{R}$.\\
For the analysis of a stochastic flow, it is convenient to divide the flow into forward flow $\varphi(t,s), s\leq t$ and the backward flow $\varphi(t,s), s\geq t$ and discuss them separately.
\end{defi}

Consider the corresponding linear ordinary differential equation:
\begin{equation}\label{21}
dX_t=AX_tdt. 
\end{equation}
Let $\Phi_j(t)=(\Phi_{1j}(t),...,\Phi_{dj}(t)) ^T $ be the solution of equation (\ref{21}) with initial value $X(0)=e_j,j=1,...,d$. Then the $d\times d$ matrix
$\Phi(t)=(\Phi_1(t),...,\Phi_d(t))=(\Phi_{ij}(t))_{d\times d}$ is the fundamental matrix of equation (\ref{21}) and $\Phi(t+s)=\Phi(t)\circ \Phi(s)$ for all $t,s\geq 0$.

Now we propose the assumptions on $A, h, \sigma$ in (\ref{15}) as follows.

$\boldsymbol{(A)}$ $A$ is cooperative, i.e., $a_{ij}\geq 0$ for all $i,j\in \{1,...,d \}$ and $i\neq j$, and $A$ is stable, i.e., all real part of its eigenvalues are negative:
\begin{equation}\label{22}
Re \ \mu\leq\lambda<0 \ \ \ \ \ \rm{for} \ \rm{all} \ \rm{eigenvalues} \ \mu \ \rm{of} \ A. 
\end{equation}
      	
$\boldsymbol{(H_1)}$ $h\in C_b^1(\mathbb{R}\times\mathbb{R}^d,\mathbb{R}_+^d)$, i.e., the function $h$ and its derivatives are both bounded, and $h$ is order-preserving in $\mathbb{R}^d$, i.e., for any $t\in \mathbb{R}$
$$h(t,x_1)\leq_{\mathbb{R}_+^d} h(t,x_2) \ \ \text{whenever} \ x_1 \leq_{\mathbb{R}_+^d} x_2,$$
or anti-order-preserving in $\mathbb{R}^d$, i.e.,
$$h(t,x_1)\geq_{\mathbb{R}_+^d} h(t,x_2) \ \ \text{whenever} \ x_1 \leq_{\mathbb{R}_+^d} x_2.$$ 
Here, $x \leq_{\mathbb{R}_+^d} y$ means that $y-x\in\mathbb{R}_+^d$ for all $x,y\in\mathbb{R}^d$.

$\boldsymbol{(H_2)}$  $L:=  \max \{\ \sup_{t\in\mathbb{R},x\in\mathbb{R}^d} |\frac{\partial h_i(t,x)}{\partial x_j}|, i,j=1,...,d\}<-\frac{\lambda}{d^2}$.   

$\boldsymbol{(\Sigma)}$   $\sup_{t\in\mathbb{R}}||\sigma(t)||:=\sup_{t\in\mathbb{R}} \max \{|\sigma_{ij}(t)| :i=1,...,d; j=1,...,m\}<\infty$, where  $\sigma_{ij}(t), i=1,...,d, j=1,...,m$ are continuous functions from $\mathbb{R}$ to $\mathbb{R}$.

If $A$ is stable, by (\cite{PM}, Chapter 2, Proposition 2.10) we have     	
\begin{equation}\label{23}
\parallel \Phi(t) \parallel:= \max \{|\Phi_{ij}(t)| :i,j=1,...,d\}\leq e^{\lambda t} , \ t\geq 0.
\end{equation}	
In this paper, we use the norm $|x|:= \max \{|x_i| :i=1,...,d\}$, $x\in\mathbb{R}^d$.

 By Theorem 4.3 in \cite{K}, there is a modification of the unique solution for equation (\ref{15}), denoted by $\varphi(t,s,\omega)x$ and the integral version of (\ref{15})
\begin{eqnarray}\label{24}
 \varphi(t,s,\omega)x
&=&\Phi(t-s)x + \Phi(t)\displaystyle \int^{t}_{s}{\Phi^{-1}(r)h(r,\varphi(r,s,\omega)x)dr}+\Phi(t)\displaystyle \int^{t}_{s}{\Phi^{-1}(r)\sigma(r) dW_r}\nonumber\\
&=&\Phi(t-s)x + \displaystyle \int^{t}_{s}{\Phi(t-r)h(r,\varphi(r,s,\omega)x)dr}+\displaystyle \int^{t}_{s}{\Phi(t-r)\sigma(r) dW_r} .
	\end{eqnarray} 
holds almost surely. In addition, it has the following properties:\\
\quad (i) For each $(t, s)\in\triangle$ and $x$, $\varphi(t,s,\cdot)x$ is ${\cal F}_s^t$-measurable.\\	
\quad (ii) For almost all $\omega$, $\varphi(t,s,\omega)x$ is continuous in $(t,s,x)$and satisfies $\lim\limits_{t\downarrow s}\varphi(t,s,x,\omega)=x$.\\
\quad (iii) For almost all $\omega$,
\begin {equation}\nonumber
\varphi(t,s,\omega)=\varphi(t,r,\omega)\circ\varphi(r,s,\omega), 
\end{equation}
for all $s\leq r\leq t$, $s,r,t\in\mathbb{R}$. 

Denote the standard $\mathbb{P}$-preserving ergodic Wiener shift by $\theta:\mathbb{R}\times\Omega\rightarrow\Omega$,
$$\theta_t\omega(\cdot):=W(t+\cdot)-W(t), \ \ \ t\in\mathbb{R}.$$

\begin{rem}
	The solution $\varphi(t,s,\omega)x$ of non-autonomous stochastic differential equation (\ref{15})  with the initial value $X(s)=x\in\mathbb{R}^d$ does not satisfy the cocycle property, but satisfies the forward stochastic flow property (\ref{12}). Furthermore, it has the periodicity (\ref{13}). In fact, for $s\leq t$, $s,t\in\mathbb{R}$, by the periodicity of $h(\cdot, x)$ we have
	\begin{equation}\nonumber
	\begin{aligned}
	\begin{split}
	& \ \  \ \ \varphi(t+T,s+T,\omega)x\\
	&=\Phi(t-s)x+ \int^{t+T}_{s+T}{\Phi(t+T-v)h(v,\varphi(v,s+T,\omega)x)dv}+\displaystyle \int^{t+T}_{s+T}{\Phi(t+T-v)\sigma(v) dW_v} \\
	&=\Phi(t-s)x+\int^{t}_{s}{\Phi(t-v)h(v+T,\varphi(v+T,s+T,\omega)x)dv}+\displaystyle \int^{t}_{s}{\Phi(t-v)\sigma(v+T) dW_{v+T}}\\
	&=\Phi(t-s)x+\int^{t}_{s}{\Phi(t-v)h(v,\varphi(v+T,s+T,\omega)x)dv}+\displaystyle \int^{t}_{s}{\Phi(t-v)\sigma(v) dW_v(\theta_{T}\omega)}
	\end{split}
	\end{aligned}
	\end{equation}
	We thus have found that the function 	  
	$$\psi(v,s,\theta_T\omega)x=\varphi(v+T,s+T,\omega)x, \ \ s\leq v $$
	satisfies
	$$\psi(t,s,\theta_T\omega)x=\Phi(t-s)x+\int^{t}_{s}{\Phi(t-v)h(v,\psi(v,s,\theta_T\omega)x)dv}+\displaystyle \int^{t}_{s}{\Phi(t-v)\sigma(v) dW_v(\theta_T\omega)}. $$	  
	By the uniqueness of the solution
	$$\varphi(t,s,\theta_T\omega)x=\psi(t,s,\theta_T\omega)x=\varphi(t+T,s+T,\omega)x.$$
	which implies that (\ref{13}) holds.
\end{rem}	  

Now we are in position to introduce the concept of an operator ${\cal K}$. Denote by $\varphi(t,-nT,\omega)x$ the solution starting from time $-nT$ with the initial value $X(-nT)=x$. Then for any $n\in\mathbb{N}_+$, $t>-nT$, 
\begin{equation}\label{25}
\varphi(t,-nT,\omega)x
=\Phi(t+nT)x + \displaystyle \int^{t}_{-nT}{\Phi(t-s)h(s,\varphi(s,-nT,\omega)x)ds}+\displaystyle \int^{t}_{-nT}{\Phi(t-s)\sigma(s) dW_s}, 
\end{equation} 
holds almost surely.

We define the input-to-state characteristic operator $\cal K$:
\begin{equation}\label{26}
[\cal{K}(u)](t,\omega)=\displaystyle \int^{t}_{-\infty}{\Phi(t-s)u(s,\omega)ds}+\displaystyle \int^{t}_{-\infty}{\Phi(t-s)\sigma(s) dW_s}, \ \ t\in\mathbb{R},\ \omega\in\Omega.
\end{equation}
where the stochastic process $u:\mathbb{R}\times\Omega\rightarrow\mathbb{R}_+^d$ is bounded.

We denote $|x|_2:=(\sum\limits_{i=1}^d {|x_i|^2})^{\frac{1}{2}}, x\in\mathbb{R}^d$ and $\parallel \Phi \parallel_2:= (\sum\limits_{i,j=1}^d {|\Phi_{i,j}|^2})^{\frac{1}{2}}$ in what follows, where $\Phi$ is a $(d\times d)$-dimensional matrix.

\begin{rem}
	It is noticed that for fixed $t\in\mathbb{R}, \omega\in\Omega$, the operator ${\cal K}$ is well defined. In fact, for any bounded stochastic process $u$, since $\parallel \Phi(t) \parallel\leq e^{\lambda t}, t\geq 0$, we have  $\parallel \Phi \parallel_2\leq d\parallel \Phi(t) \parallel\leq d e^{\lambda t}$, and it is not hard to prove that 
	$$ \displaystyle \int^{t}_{-\infty}|{\Phi(t-s)u(s,\omega)|_2ds} < \infty $$
	which implies that $\lim_{r\rightarrow\infty}\int^{t}_{-r}{\Phi(t-s)u(s,\omega)ds}$ exists for all $t\in\mathbb{R}, \omega\in\Omega$. For any $r_1>r_2>0$,
	\begin{equation}\nonumber
	\begin{aligned}
	\begin{split}
	&\ \ \ \mathbb{E}\Big|\displaystyle \int^{t}_{-r_1}{\Phi(-s)\sigma(s) dW_s}-\displaystyle \int^{t}_{-r_2}{\Phi(-s)\sigma(s) dW_s}\Big|_2^2\\
	&=\mathbb{E}\displaystyle \int^{-r_2}_{-r_1}{||\Phi(-s)\sigma(s)||_2^2 ds} \\
	&\leq  d m\cdot \sup\limits_{t\in\mathbb{R}}||\sigma(t)||\cdot\displaystyle \int^{-r_2}_{-r_1}{||\Phi(-s)||_2^2 ds}
	\end{split}
	\end{aligned}
	\end{equation}
	which together with assumption $\boldsymbol{(A)}$ and $\boldsymbol{(\Sigma)}$ shows that $\int^{t}_{-r}{\Phi(-s)\sigma(s) dW_s}$ converges in $L^2$ as $r\rightarrow \infty$. Furthermore, $\int^{t}_{-r}{\Phi(-s)\sigma(s) dW_s}$ is a continuous martingale. Hence, it follows from (\cite{KS}, Problem 3.20 in Chapter 1) that $\int^{t}_{-r}{\Phi(-s)\sigma(s) dW_s}$ converges $\mathbb{P}$-a.s. to an integrable random variable $X_{\infty}:=\int^{t}_{-\infty}{\Phi(-s)\sigma(s) dW_s}$ as $r\rightarrow\infty$.
\end{rem}     	

\section{Measurability and asymptotic behaviour}
\quad In this section, we give some lemmas to describe the dynamics of the pull-back trajectory which will be used in the proof of our main result.


\begin{lem}
	For each $n\in\mathbb{N}^+$, $x\in\mathbb{R}^d$, let
	$$a_n^h(t, \omega)=\inf{\{h(t, \varphi(t,-mT,\omega)x):m\geq n,m\in\mathbb{N}_+\}},\ \ t\in\mathbb{R}, \omega\in\Omega,$$
	and
	$$b_n^h(t, \omega)=\ \sup{\{h(t, \varphi(t,-mT, \omega)x):m\geq n, m\in\mathbb{N}_+\}},\ \ t\in\mathbb{R}, \omega\in\Omega,$$
	where $\inf$ and $\sup$ mean the greatest lower bound and the least upper bound, respectively.
	Then $a_n^h(t, \omega)$ and $b_n^h(t, \omega)$ are $({\cal B}(\mathbb{R})\otimes{\cal F},{\cal B}(\mathbb{R}^d))$-measurable stochastic processes. 
	
\end{lem}

\noindent {\bf Proof.}\ \ First, we show that $a_n^h(t, \omega)$ and $b_n^h(t,\omega)$ are well defined. It is clear that 
$$D_n^h(t, \omega):={\{h(t, \varphi(t,-mT,\omega)x):m\geq n,m\in\mathbb{N}_+\}}$$ 
is a bounded set for fixed $n\in\mathbb{N}_+, t\in\mathbb{R}, \omega\in\Omega$, which implies that $D_n^h(t, \omega)$ is order-bounded. Since $\mathbb{R}_+^d$ is strongly minihedral (\cite{C}, Definition 3.1.7), $a_n^h(t, \omega)$ and $b_n^h(t,\omega)$ exist. Since $\inf {D_n^h(t, \omega)}=-\sup \{-D_n^h(t, \omega)\}$, it is sufficient to consider $b_n^h(t, \omega)$ only. Let
$$\beta_{n,M}^h(t,\omega)=\sup\Big\{h(t, \varphi(t,-nT,\omega)x), h(t, \varphi(t,-(n+1)T,\omega)x), \cdots, h(t, \varphi(t,-(n+M)T,\omega)x) \Big\},$$
by the continuity of $h$, the measurability of $\varphi$ and Corollary 3.1.1(ii) in \cite{C}, $\beta_{n,M}^h(t,\omega)$ is $({\cal B}(\mathbb{R})\otimes{\cal F},{\cal B}(\mathbb{R}^d))$-measurable for every $M=1,2,\cdots$. It is also clear that
$$\beta_{n,1}^h(t,\omega)\leq \beta_{n,2}^h(t,\omega)\leq \cdots\leq \beta_{n,M}^h(t,\omega)\leq\cdots .$$
Moreover, by the boundedness of $h$ in $\mathbb{R}_+^d$, $b_n^h(t, \omega)=\lim\limits_{{M} \to {\infty}}\beta_{n,M}^h(t,\omega)$ is $({\cal B}(\mathbb{R})\otimes{\cal F},{\cal B}(\mathbb{R}^d))$-measurable stochastic process. 
\hfill\fbox\\

\begin{lem}
	Assume that conditions $\boldsymbol{(A)}$, $\boldsymbol{(H_1)}$ and $\boldsymbol{(\Sigma)}$ hold. Let $\varphi(t,s,\omega)x$ be a solution of stochastic system (\ref{15}) with initial value $X(s)=x\in\mathbb{R}^d$. Then we have
	
	\begin{equation}\label{31}
	{\cal K}(\underline{\lim}h(\cdot, \varphi))\leq  \underline{\lim}\varphi\leq  \overline{\lim}\varphi\leq {\cal K}( \overline{\lim}h(\cdot, \varphi)),\ \ \mathbb{P}-a.s.,
	\end{equation}  
	where 
	$$[ \underline{\lim}\varphi](t,\omega)=\lim \limits_{{n} \to {\infty}} \inf\{\varphi(t, -mT, \omega)x:m\geq n, m\in\mathbb{N}_+\}, \ \ x\in\mathbb{R}^d, t\in\mathbb{R}, \omega\in\Omega.$$
	
	$$[ \overline{\lim}\varphi](t,\omega)=\lim \limits_{{n} \to {\infty}}\sup\{\varphi(t, -mT, \omega)x:m\geq n, m\in\mathbb{N}_+\}, \ \ x\in\mathbb{R}^d, t\in\mathbb{R}, \omega\in\Omega.$$
	
	$$[ \underline{\lim}h(\cdot, \varphi)](t,\omega)=\lim \limits_{{n} \to {\infty}} a_n^h(t,\omega), \ \ x\in\mathbb{R}^d, t\in\mathbb{R}, \omega\in\Omega.$$
	and 
	$$[ \overline{\lim}h(\cdot, \varphi)](t,\omega)=\lim \limits_{{n} \to {\infty}}b_n^h(t,\omega), \ \ x\in\mathbb{R}^d, t\in\mathbb{R}, \omega\in\Omega.$$
\end{lem}
where $a_n^h$ and $b_n^h$ are as defined in Lemma 3.1.

\noindent {\bf Proof.}\ \ Here, we only prove the first inequality for the sake of convenience and the rest of the inequalities can be proved analogously. Similar to Lemma 3.1, we can easily get that $ \underline{\lim}\varphi$ and $ \underline{\lim}h(\cdot, \varphi)$ exist, which are also two $({\cal B}(\mathbb{R})\otimes{\cal F},{\cal B}(\mathbb{R}^d))$-measurable stochastic processes. Then by (\ref{26}) and the Fubini Theorem, ${\cal K}( \underline{\lim}h(\cdot, \varphi))$ is well defined and $({\cal B}(\mathbb{R})\otimes{\cal F},{\cal B}(\mathbb{R}^d))$-measurable. By the definition of the $ \underline{\lim}h(\cdot, \varphi)$, ${\cal K}$ and Lebesgue's dominated convergence theorem, we have ${\cal K}( \underline{\lim}h(\cdot, \varphi))=\lim\limits_{{n} \to {\infty}}{\cal K}(a_n^h)$. For fixed $n\in\mathbb{N}_+$, it is enough to prove that 
\begin{equation}\nonumber
\begin{aligned}
\begin{split}
&\ \ \ \ [{\cal K}(a_n^h)](t, \omega)\\
& = \displaystyle \int^{t}_{-\infty}{\Phi(t-s)\inf\{h(s,\varphi(s, -mT, \omega)x):m\geq n, m\in\mathbb{N}_+\}ds}+\displaystyle \int^{t}_{-\infty}{\Phi(t-s)\sigma(s) dW_s}\\
&=\lim \limits_{{\tilde{m}} \to {\infty}\atop{\tilde{m}}\geq n} \Big\{\Phi(t+{\tilde{m}}T)x +\displaystyle \int^{t}_{n T-{\tilde{m}T}}{\Phi(t-s)\inf\{h(s,\varphi(s, -mT, \omega)x):m\geq n, m\in\mathbb{N}_+\}ds}\\
& \ \ \ \ +\displaystyle \int^{t}_{-{\tilde{m}T}}{\Phi(t-s)\sigma(s) dW_s}\Big\}\\
&=\lim \limits_{\tilde{n} \to {\infty}\atop \tilde{n}\geq n}\inf \Big\{\Phi(t+{\tilde{m}}T)x +\displaystyle \int^{t}_{n T-{\tilde{m}T}}{\Phi(t-s)\inf\{h(s,\varphi(s, -mT, \omega)x):m\geq n, m\in\mathbb{N}_+\}ds}\\
& \ \ \ \ +\displaystyle \int^{t}_{-{\tilde{m}T}}{\Phi(t-s)\sigma(s) dW_s}:{\tilde{m}}\geq \tilde{n}, {\tilde{m}}\in\mathbb{N}_+\Big\}\\
&\leq \lim \limits_{\tilde{n} \to {\infty}\atop \tilde{n}\geq n}\inf \Big\{\Phi(t+{\tilde{m}}T)x +\displaystyle \int^{t}_{n T-{\tilde{m}T}}{\Phi(t-s)h(s,\varphi(s, -{\tilde{m}}T, \omega)x)ds}\\
& \ \ \ \ +\displaystyle \int^{t}_{-{\tilde{m}T}}{\Phi(t-s)\sigma(s) dW_s}:{\tilde{m}}\geq \tilde{n},{\tilde{m}}\in\mathbb{N}_+\Big\}\\    
&\leq \lim \limits_{\tilde{n} \to {\infty}}\inf \Big\{\Phi(t+{\tilde{m}}T)x +\displaystyle \int^{t}_{-{\tilde{m}T}}{\Phi(t-s)h(s,\varphi(s, -{\tilde{m}}T, \omega)x)ds}\\
& \ \ \ \ +\displaystyle \int^{t}_{-{\tilde{m}T}}{\Phi(t-s)\sigma(s) dW_s}:{\tilde{m}}\geq \tilde{n},{\tilde{m}}\in\mathbb{N}_+\Big\}\\ 
&= [ \underline{\lim}\varphi](t,\omega),
\end{split}
\end{aligned}
\end{equation}
where the third equality has used Lemma A.2 in \cite{MS}, while the second-to-last inequality has applied the positivity of $\Phi(t)$ and $h$.
\hfill\fbox\\   
\begin{lem}
	Assume that conditions $\boldsymbol{(A)}$, $\boldsymbol{(H_1)}$ and $\boldsymbol{(\Sigma)}$ hold. Let $\varphi(t,s,\omega)x$ be a solution of stochastic system (\ref{15}) with initial value $X(s)=x\in\mathbb{R}^d$. Then we have the following:\\
	$(i)$ if $h$ is order-preserving in $\mathbb{R}^d$, then
	\begin{equation}\label{32}
	h(\cdot, \underline{\lim}\varphi)\leq  \underline{\lim}h(\cdot,\varphi)\leq  \overline{\lim}h(\cdot,\varphi)\leq h(\cdot,  \overline{\lim}\varphi),\ \ \mathbb{P}-a.s.,
	\end{equation}  
	$(ii)$ if $h$ is anti-order-preserving in $\mathbb{R}^d$, then
	\begin{equation}\label{33}
	h(\cdot,  \overline{\lim}\varphi)\leq  \underline{\lim}h(\cdot,\varphi)\leq  \overline{\lim}h(\cdot,\varphi)\leq  h(\cdot, \underline{\lim}\varphi),\ \ \mathbb{P}-a.s.,
	\end{equation} 
\end{lem}
\noindent {\bf Proof.}\ \ Indeed, the proof of the first inequality in (\ref{32}) is adequate and the rest of the results of this lemma can be obtained analogously. Observe that $h$ is order-preserving in $\mathbb{R}^d$; then for fixed $n\in \mathbb{N}_+, m\geq n, m\in\mathbb{N}_+$, we have
$$h(t,\inf\{\varphi(t, -kT, \omega)x:k\geq n, k\in\mathbb{N}_+\})\leq h(t,\varphi(t, -mT, \omega)x)$$
and
\begin{equation}\label{34}
h(t,\inf\{\varphi(t, -mT, \omega)x:m\geq n,m\in\mathbb{N}_+\})\leq \inf \{h(t,\varphi(t, -mT, \omega)x):m\geq n,m\in\mathbb{N}_+\}. 
\end{equation} 
Let $n\rightarrow\infty$ in (\ref{34}). By the continuity of $h$, we have
\begin{equation}\nonumber
\begin{aligned}
\begin{split}
[h(\cdot, \underline{\lim}\varphi)](t,\omega)
&= h(t,\lim \limits_{{n} \to {\infty}}\inf\{\varphi(t, -mT, \omega)x:m\geq n,m\in\mathbb{N}_+\})\\
&= \lim \limits_{{n} \to {\infty}}h(t,\inf\{\varphi(t, -mT, \omega)x:m\geq n,m\in\mathbb{N}_+\})\\
&\leq  \lim \limits_{{n} \to {\infty}} \inf\{h(t,\varphi(t, -mT, \omega)x):m\geq n,m\in\mathbb{N}_+\}\\
&=[ \underline{\lim}h(\cdot, \varphi)](t,\omega).
\end{split}
\end{aligned}
\end{equation}   
\hfill\fbox\\  
\begin{lem}
	Assume that conditions $\boldsymbol{(A)}$, $\boldsymbol{(H_1)}$ and $\boldsymbol{(\Sigma)}$ hold. Let $\varphi(t,s,\omega)x$ be a solution of stochastic system (\ref{15}) with initial value $X(s)=x\in\mathbb{R}^d$. Then we have
	\begin{equation}\label{35}
	{\cal K}(a_n^h)\leq \underline{\lim}\varphi\leq \overline{\lim}\varphi\leq{\cal K}(b_n^h) \ \ \mathbb{P}-a.s.,n\in\mathbb{N}_+,
	\end{equation}    
	where $a_n^h$ and $b_n^h$ are as defined in Lemma 3.1. Furthermore, define the gain operator
	$${\cal K}^h(u)(t, \omega)=h(t,[{\cal K}(u)](t, \omega)).$$ Then we have the following:\\
	$(i)$ if $h$ is order-preserving in $\mathbb{R}^d$, then for fixed $n\in\mathbb{N}_+$,
	\begin{equation}\label{36}
	({\cal K}^h)^k(a_n^h)\leq  \underline{\lim}h(\cdot,\varphi)\leq  \overline{\lim}h(\cdot,\varphi)\leq ({\cal K}^h)^k(b_n^h),\ \ \mathbb{P}-a.s.,\ \ k\in\mathbb{N}_+.
	\end{equation}  
	$(ii)$ if $h$ is anti-order-preserving in $\mathbb{R}^d$, then for fixed $n\in\mathbb{N}_+$,
	\begin{equation}\label{37}
	({\cal K}^h)^{2k}(a_n^h)\leq  \underline{\lim}h(\cdot,\varphi)\leq  \overline{\lim}h(\cdot,\varphi)\leq ({\cal K}^h)^{2k}(b_n^h),\ \ \mathbb{P}-a.s., \ \ k\in\mathbb{N}_+.
	\end{equation} 
\end{lem}

\noindent {\bf Proof.}\ \  By the definition of $a_n^h$ and $b_n^h$, it is evident that 
$$a_n^h\leq  \underline{\lim}h(\cdot,\varphi)\leq  \overline{\lim}h(\cdot,\varphi)\leq b_n^h,\ \ n\in\mathbb{N}_+.$$
By the positivity of $\Phi$, ${\cal K}(u)$ is monotone with respect to $u$, and consequently   
$${\cal K}(a_n^h)\leq {\cal K}( \underline{\lim}h(\cdot,\varphi))\leq {\cal K}(  \overline{\lim}h(\cdot,\varphi))\leq {\cal K}(b_n^h), \ \ \mathbb{P}-a.s., \ \  n\in\mathbb{N}_+.$$
by (\ref{31}), we have
$${\cal K}(a_n^h)\leq  \underline{\lim}\varphi\leq  \overline{\lim}\varphi\leq {\cal K}(b_n^h), \ \ \mathbb{P}-a.s., \ \  n\in\mathbb{N}_+.$$
which implies that (\ref{35}) holds. 

In what follows, we claim that (\ref{36}) and (\ref{37}) hold.

If $h$ is order-preserving in $\mathbb{R}^d$, then it deduces that $h$ preserves the inequalities in (\ref{35}):
$${\cal K}^h(a_n^h)\leq h(\cdot, \underline{\lim}\varphi)\leq h(\cdot, \overline{\lim}\varphi)\leq {\cal K}^h(b_n^h), \ \ \mathbb{P}-a.s., \ \  n\in\mathbb{N}_+.$$ 
which together with (\ref{32}) implies
$${\cal K}^h(a_n^h)\leq  \underline{\lim}h(\cdot,\varphi)\leq  \overline{\lim}h(\cdot,\varphi)\leq {\cal K}^h(b_n^h), \ \ \mathbb{P}-a.s., \ \  n\in\mathbb{N}_+.$$  
This proves that (\ref{36}) is true for $k=1$.

Next we assume that, for some $k\in\mathbb{N}$, we have obtained
$$({\cal K}^h)^k(a_n^h)\leq  \underline{\lim}h(\cdot,\varphi)\leq  \overline{\lim}h(\cdot,\varphi)\leq ({\cal K}^h)^k(b_n^h), \ \ \mathbb{P}-a.s., \ \  n\in\mathbb{N}_+.$$
From the monotonicity of ${\cal K}$ and (\ref{31}), we have:
\begin{equation}\nonumber
\begin{aligned}
\begin{split}
{\cal K}[({\cal K}^h)^k(a_n^h)]
&\leq {\cal K}( \underline{\lim}h(\cdot,\varphi))\leq  \underline{\lim}\varphi\\
&\leq  \overline{\lim}\varphi\leq {\cal K}( \overline{\lim}h(\cdot,\varphi))\leq  {\cal K}[({\cal K}^h)^k(b_n^h)].
\end{split}
\end{aligned}
\end{equation} 
By the monotonicity of $h$ in $\mathbb{R}^d$ and (\ref{32}), we get that
$$({\cal K}^h)^{k+1}(a_n^h)\leq  \underline{\lim}h(\cdot,\varphi)\leq  \overline{\lim}h(\cdot,\varphi)\leq ({\cal K}^h)^{k+1}(b_n^h), \ \ \mathbb{P}-a.s., \ \  n\in\mathbb{N}_+.$$ 
Therefore, we conclude that (\ref{36}) holds by mathematical induction.

If $h$ is anti-order-preserving in $\mathbb{R}^d$, similar to $h$ is order-preserving in $\mathbb{R}^d$; we deduce that
$${\cal K}^h(b_n^h)\leq h(\cdot, \overline{\lim}\varphi)\leq h(\cdot, \underline{\lim}\varphi)\leq {\cal K}^h(a_n^h), \ \ \mathbb{P}-a.s., \ \  n\in\mathbb{N}_+.$$
by (\ref{33}), we have
$${\cal K}^h(b_n^h)\leq \underline{\lim}h(\cdot,\varphi)\leq  \overline{\lim}h(\cdot,\varphi)\leq {\cal K}^h(a_n^h), \ \ \mathbb{P}-a.s., \ \  n\in\mathbb{N}_+.$$
Combining the monotonicity of ${\cal K}$ and (\ref{31}), it shows that
$${\cal K}[{\cal K}^h(b_n^h)]\leq \underline{\lim}\varphi\leq  \overline{\lim}\varphi\leq {\cal K}[{\cal K}^h(a_n^h)], \ \ \mathbb{P}-a.s., \ \  n\in\mathbb{N}_+,$$
which together with the anti-monotonicity of $h$ in $\mathbb{R}^d$ and (\ref{33}) implies
$$({\cal K}^h)^2(a_n^h)\leq  \underline{\lim}h(\cdot,\varphi)\leq  \overline{\lim}h(\cdot,\varphi)\leq ({\cal K}^h)^2(b_n^h), \ \ \mathbb{P}-a.s., \ \  n\in\mathbb{N}_+.$$ 
The rest of the proof of (\ref{37}) can be obtained analogously to $h$ is order-preserving in $\mathbb{R}^d$ by the mathematical induction.
\hfill\fbox\\

\section{Main results} 
\quad In this section, we state our main result on the existence of random periodic solution of nonlinear stochastic system (\ref{15}) and present its proof. We begin with a lemma.  

Let ${\cal M}_{{\cal B}(\mathbb{R})\otimes{\cal F}}^{b, T}(\mathbb{R}\times\Omega;[0,N])$ be the space of ${\cal B}(\mathbb{R})\otimes{\cal F}$-measurable functions $f:\mathbb{R}\times\Omega\rightarrow [0,N]$ with $f(t+T, \omega)\leq f(t, \theta_T\omega), t\in\mathbb{R}, \omega\in\Omega$ or $f(t+T, \omega)\geq f(t, \theta_T\omega), t\in\mathbb{R}, \omega\in\Omega$, where $N=(N_1,...,N_d), N_i= \sup_{t\in\mathbb{R}, x\in\mathbb{R}^d}|h_i(t,x)|, i=1,...,d$. We introduce a metric on ${\cal M}_{{\cal B}(\mathbb{R})\otimes{\cal F}}^{b, T}(\mathbb{R}\times\Omega;[0,N])$ as follows:
$$\varrho(f_1, f_2):=|f_1-f_2|_\infty=\mathop{ \sup}\limits_{t\in\mathbb{R}, \omega\in\Omega}{|f_1(t,\omega)-f_2(t,\omega)|},  \ \ \text{for} \ \text{all}\ f_1, f_2\in{\cal M}_{{\cal B}(\mathbb{R})\otimes{\cal F}}^{b, T}.$$

\begin{lem}
	Assume that conditions $\boldsymbol{(A)}$, $\boldsymbol{(H_1)}$, $\boldsymbol{(H_2)}$ and $\boldsymbol{(\Sigma)}$ hold. Then  $({\cal M}_{{\cal B}(\mathbb{R})\otimes{\cal F}}^{b, T}, \varrho)$ is a complete metric space and the gain operator ${\cal K}^h: {\cal M}_{{\cal B}(\mathbb{R})\otimes{\cal F}}^{b, T}\rightarrow {\cal M}_{{\cal B}(\mathbb{R})\otimes{\cal F}}^{b, T}$, $u\mapsto{\cal K}^h(u) $ is a contractive mapping, where ${\cal K}^h(u)(t, \omega)=h(t,[{\cal K}(u)](t, \omega))$ and the definition of the input-to-state characteristic operator ${\cal K}$ can be chosen as an $\mathbb{R}^d$-value version for all $t\in\mathbb{R}, \omega\in\Omega$. 
\end{lem}
\noindent {\bf Proof.}\ \ It is clear that $({\cal M}_{{\cal B}(\mathbb{R})\otimes{\cal F}}^{b, T}, \varrho)$ is a metric space. We show that the metric space  ${\cal M}_{{\cal B}(\mathbb{R})\otimes{\cal F}}^{b, T}$ is complete with respect to $\varrho$.
To prove this, we choose a Cauchy sequence $\{f_n, n\in\mathbb{N}\}$ in $({\cal M}_{{\cal B}(\mathbb{R})\otimes{\cal F}}^{b, T}, \varrho)$; we denote a function $f$ as follows:
$$f(t, \omega):=\lim \limits_{{n} \to {\infty}}{f_n(t, \omega)}\in [0,N] \ \ \text{for} \ \text{all} \ t\in\mathbb{R},\omega\in\Omega,$$
which holds based on the fact that $\{f_n(t,\omega), n\in \mathbb{N}\}$ is a Cauchy sequence in $\mathbb{R}^d$ for fixed $t\in\mathbb{R}, \omega\in\Omega$. It is noticed that the limit of a family of ${\cal B}(\mathbb{R})\otimes{\cal F}$-measurable functions is an ${\cal B}(\mathbb{R})\otimes{\cal F}$-measurable function. Without loss of generality, we assume that there is a subsequence $\{f_{n_k}\}$ of $\{f_n\}$ such that
$$f_{n_k}(t+T, \omega)\leq f_{n_k}(t,\theta_T\omega),\ \  t\in\mathbb{R},\  \omega\in\Omega.$$
Let $n_k\rightarrow \infty$ in the above inequality, then $f(t+T,\omega)\leq f(t,\theta_T\omega), t\in\mathbb{R}, \omega\in\Omega$. In what follows, we will prove that $|f-f_n|_{\infty}\rightarrow 0$ as $n\rightarrow \infty$. It is noticed that $\{f_n(t, \omega), n\in N \}$ is a Cauchy sequence, we know that for any $\varepsilon>0$, there exists an $N_0=N_0(\varepsilon)\in \mathbb{N}$ such that for $n, m\geq N_0$,
$$\mathop{ \sup}\limits_{t\in\mathbb{R}, \omega\in\Omega}{|f_m(t,\omega)-f_n(t,\omega)|}<\varepsilon.$$
Let $m\rightarrow \infty$; then
$$\mathop{ \sup}\limits_{t\in\mathbb{R}, \omega\in\Omega}{|f(t,\omega)-f_n(t,\omega)|}\leq\varepsilon \ \ \text{for} \ \text{all} \ n\geq N_0,$$
which implies that $|f-f_n|_{\infty}\rightarrow 0$ as $n\rightarrow \infty$. Thus $({\cal M}_{{\cal B}(\mathbb{R})\otimes{\cal F}}^{b, T}, \varrho)$ is a complete metric space.

Next we claim that ${\cal K}^h:{\cal M}_{{\cal B}(\mathbb{R})\otimes{\cal F}}^{b, T}\rightarrow{\cal M}_{{\cal B}(\mathbb{R})\otimes{\cal F}}^{b, T}$ is a contractive mapping. First, we should show that  ${\cal K}^h:{\cal M}_{{\cal B}(\mathbb{R})\otimes{\cal F}}^{b, T}\rightarrow{\cal M}_{{\cal B}(\mathbb{R})\otimes{\cal F}}^{b, T}$ is well defined. From $\boldsymbol{(H_1)}$, it follows that $h:\mathbb{R}\times\mathbb{R}^d\rightarrow [0, N]$ and $h(t+T, \omega)= h(t, \omega), t\in\mathbb{R}, \omega\in\Omega$. For any $f\in ({\cal M}_{{\cal B}(\mathbb{R})\otimes{\cal F}}^{b, T}, \varrho)$, by the definition of ${\cal K}$, the measurability of $\theta$, and the Fubini theorem, it is evident that ${\cal K}(f)$ is a ${\cal B}(\mathbb{R})\otimes{\cal F}$-measurable function. Without loss of generality,  we assume that $f(t+T,\omega)\leq f(t,\theta_T\omega), t\in\mathbb{R}, \omega\in\Omega$, then
\begin{equation}\nonumber
\begin{aligned}
\begin{split}
[{\cal K}(f)](s,\theta_T\omega)
& = \displaystyle \int^{s}_{-\infty}{\Phi(s-r)f(r,\theta_T\omega))dr}+\displaystyle \int^{s}_{-\infty}{\Phi(s-r)\sigma(r) dW_r({\theta_T}\omega)}\\
&=\displaystyle \int^{s+T}_{-\infty}{\Phi(s+T-r)f(r-T,\theta_T\omega)dr}+\displaystyle \int^{s+T}_{-\infty}{\Phi(s+T-r)\sigma(r-T) dW_r}\\
&\geq \displaystyle \int^{s+T}_{-\infty}{\Phi(s+T-r)f(r,\omega)dr}+\displaystyle \int^{s+T}_{-\infty}{\Phi(s+T-r)\sigma(r) dW_r}\\
&= [{\cal K}(f)](s+T,\omega) \ \ \  \mathbb{P}-a.s.
\end{split}
\end{aligned}
\end{equation}
If $h$ is order-preserving in $\mathbb{R}^d$, we have
\begin{equation}\nonumber
\begin{aligned}
\begin{split}
{\cal K}^h(f)(t+T, \omega)
& = h(t+T,{\cal K}(f)(t+T, \omega))\\
&\leq h(t,{\cal K}(f)(t, \theta_T\omega))\\
&= {\cal K}^h(f)(t, \theta_T\omega)
\end{split}
\end{aligned}
\end{equation}
Otherwise, if $h$ is anti-order-preserving in $\mathbb{R}^d$, we have ${\cal K}^h(f)(t+T, \omega)\geq {\cal K}^h(f)(t, \theta_T\omega)$. Which yields ${\cal K}^h:{\cal M}_{{\cal B}(\mathbb{R})\otimes{\cal F}}^{b, T}\rightarrow{\cal M}_{{\cal B}(\mathbb{R})\otimes{\cal F}}^{b, T}$. Finally, we prove that ${\cal K}^h$ is a contractive mapping. By $\boldsymbol{(H_1)}$ and $\boldsymbol{(H_2)}$, we have
$$ \mathop{ \sup}\limits_{t\in\mathbb{R}, x\in\mathbb{R}^d}{||D_x(h(t,x))||}=\mathop{ \sup}\limits_{t\in\mathbb{R}, x\in\mathbb{R}^d}{\Big|\Big|\Big(\frac{\partial h_i(t,x)}{\partial x_j}\Big)_{d\times d}\Big|\Big|}\leq L,$$
Let $f_1, f_2$ be two elements in $({\cal M}_{{\cal B}(\mathbb{R})\otimes{\cal F}}^{b, T}, \varrho)$. By the fact that $|\Phi(x)|\leq d||\Phi||\cdot|x|$ for all $x\in\mathbb{R}^d$ and $\Phi\in\mathbb{R}^{d\times d}$, we get
\begin{equation}\nonumber
\begin{aligned}
\begin{split}
|{\cal K}^h(f_1)-{\cal K}^h(f_2)|_{\infty}
& = |h(\cdot,{\cal K}(f_1))-h(\cdot,{\cal K}(f_2))|_{\infty}\\
& = |D_xh(\cdot,[{\cal K}(f_2)+\mu({\cal K}(f_1)-{\cal K}(f_2))])\cdot[{\cal K}(f_1)-{\cal K}(f_2)]|_{\infty}\\
&\leq  d\mathop{ \sup}\limits_{t\in\mathbb{R}, x\in\mathbb{R}^d}{||D_x(h(t,x))||}\cdot|{\cal K}(f_1)-{\cal K}(f_2)|_{\infty}\\
&\leq Ld\Big|\displaystyle \int^{t}_{-\infty}{\Phi(t-s)f_1(s,\omega)}ds-\displaystyle \int^{t}_{-\infty}{\Phi(t-s)f_2(s,\omega))}ds\Big|_{\infty}\\
&\leq Ld^2\displaystyle \int^{t}_{-\infty}{||\Phi(t-s)||\cdot|f_1-f_2|_{\infty}}ds\\
&\leq Ld^2|f_1-f_2|_{\infty}\displaystyle \int^{t}_{-\infty}{e^{\lambda(t-s)}}ds\\
&=-\frac{Ld^2}{\lambda}|f_1-f_2|_{\infty}
\end{split}
\end{aligned}
\end{equation}
where $\mu\in(0,1)$,$-\frac{Ld^2}{\lambda}<1$.

\hfill\fbox\\

\begin{thm}
	Assume that conditions $\boldsymbol{(A)}$, $\boldsymbol{(H_1)}$, $\boldsymbol{(H_2)}$ and $\boldsymbol{(\Sigma)}$ hold. Then the gain operator 
	$${\cal K}^h: {\cal M}_{{\cal B}(\mathbb{R})\otimes{\cal F}}^{b, T}\rightarrow {\cal M}_{{\cal B}(\mathbb{R})\otimes{\cal F}}^{b, T}, \ \ u\mapsto{\cal K}^h(u) $$
	possesses a unique nonnegative fixed point $u\in {\cal M}_{{\cal B}(\mathbb{R})\otimes{\cal F}}^{b, T}(\mathbb{R}\times\Omega;[0,N])$ such that for all $t\in \mathbb{R}, x\in\mathbb{R}^d$,
	\begin{equation}\label{41}
	\lim \limits_{{n} \to {\infty}}{\varphi(t, -nT, \omega)x}=[{\cal K}(u)](t, \omega),\ \ \mathbb{P}-a.s.
	\end{equation}
	Moreover, $\varphi(t,s,\omega)[{\cal K}(u)](s,\omega)=[{\cal K}(u)](t,\omega)$, $[{\cal K}(u)](s+T,\omega)=[{\cal K}(u)](s,\theta_T\omega)$, $\mathbb{P}-a.s.$, $(t, s)\in\triangle$; i.e., the image $[{\cal K}(u)](\cdot,\cdot)$ at the fixed point $u$ for the input-to-state characteristic operator ${\cal K}$ is a random periodic solution. 
\end{thm}
\noindent {\bf Proof.}\ \
In view of Lemma 3.4, regardless of the monotonicity or anti-monotonicity for h, for fixed $n\in\mathbb{N}_+$, we have
\begin{equation}\label{42}
({\cal K}^h)^{2k}(a_n^h)\leq \underline\lim h(\cdot, \varphi)\leq\overline\lim h(\cdot, \varphi)\leq ({\cal K}^h)^{2k}(b_n^h) \ \ \mathbb{P}-a.s., k\in\mathbb{N}
\end{equation} 
where $a_n^h$ and $b_n^h$ are as defined in Lemma 3.1. By Lemma 3.1, $a_n^h$ and $b_n^h$ are bounded ${\cal B}(\mathbb{R})\otimes{\cal F}$-measurable functions. By (\ref{13}) and the definition of $a_n^h$ in Lemma 3.1, we have
\begin{equation}\nonumber
\begin{aligned}
\begin{split}
a_n^h(t+T, \omega)
& = \inf{\{h(t+T, \varphi(t+T,-mT,\omega)x):m\geq n,m\in\mathbb{N}_+\}}\\
& = \inf{\{h(t, \varphi(t,-(m+1)T,\theta_T\omega)x):m\geq n,m\in\mathbb{N}_+\}}\\
&\geq \inf{\{h(t, \varphi(t,-mT,\theta_T\omega)x):m\geq n,m\in\mathbb{N}_+\}}\\
&= a_n^h(t, \theta_T\omega).
\end{split}
\end{aligned}
\end{equation}
Similarly, we have $b_n^h(t+T, \omega)\leq b_n^h(t, \theta_T\omega)$. So $a_n^h$ and $b_n^h$ are both in $({\cal M}_{{\cal B}(\mathbb{R})\otimes{\cal F}}^{b, T}, \varrho)$. Since ${\cal K}^h$ is a contractive mapping on the complete metric space $({\cal M}_{{\cal B}(\mathbb{R})\otimes{\cal F}}^{b, T}, \rho)$, by the Banach fixed point theorem (\cite{Y}), there exists a unique nonnegative stochastic process $u: \mathbb{R}\times\Omega\rightarrow [0,N]$ for ${\cal K}^h$ such that
$$[{\cal K}^h(u)](t, \omega)=u(t, \omega) \ \ \ \text{for} \ \text{all} \ t\in\mathbb{R}, \omega\in\Omega.$$
then
\begin{equation}\label{43}
\lim \limits_{{k} \to {\infty}}{[({\cal K}^h)^{2k}(a_n^h)](t, \omega)}=u(t, \omega)= \lim \limits_{{k} \to {\infty}}{[({\cal K}^h)^{2k}(b_n^h)](t, \omega)} \ \ \ \rm{for} \ \rm{all} \ t\in\mathbb{R}, \omega\in\Omega.
\end{equation}

It is noticed that $b_n^h(t+T, \omega)\leq b_n^h(t, \theta_T\omega)$, $h(t+T, x)=h(t,x)$ and $({\cal K}^h)^{2k}$ are monotone, then
$$[({\cal K}^h)^{2k}(b_n^h)](t+T, \omega)\leq [({\cal K}^h)^{2k}(b_n^h)](t, \theta_T\omega)$$
Let $k\rightarrow\infty$ in the above inequality, then $u(t+T,\omega)\leq u(t,\theta_T\omega)$. Similarly, for $a_n^h(t+T, \omega)\geq a_n^h(t, \theta_T\omega)$ we have $u(t+T,\omega)\geq u(t,\theta_T\omega)$. So for all $t\in\mathbb{R}, \omega\in\Omega$
\begin{equation}\label{44}
u(t+T,\omega)= u(t,\theta_T\omega)
\end{equation}    
Combining (\ref{42}) and (\ref{43}), we have
$$[ \underline\lim h(\cdot, \varphi)](t, \omega)=[ \overline\lim h(\cdot, \varphi)](t, \omega)=u(t, \omega) \ \ \mathbb{P}-a.s.$$
which together with (\ref{31}) implies that 
$$[ \underline\lim\varphi](t, \omega)=[ \overline\lim \varphi](t, \omega)=[{\cal K}(u)](t, \omega) \ \ \mathbb{P}-a.s.$$ 
In order to prove (\ref{41}), it remains to show that
\begin{equation}\label{45}
[ \underline\lim\varphi](t, \omega)=[ \overline\lim \varphi](t, \omega)=\lim \limits_{{n} \to {\infty}}\varphi(t,-nT, \omega)x \ \ \mathbb{P}-a.s., x\in\mathbb{R}^d.
\end{equation}
By the definition of infimum and supremum, it is clear that
\begin{equation}\nonumber
\begin{aligned}
\begin{split}
\inf\{\varphi(t, -mT, \omega)x:m\geq n \}
&\leq {\varphi(t,-nT, \omega)x}\\
&\leq \ \sup\{\varphi(t, -mT, \omega)x:m\geq n \} \ \ \mathbb{P}-a.s., x\in\mathbb{R}^d.
\end{split}
\end{aligned}
\end{equation}
Let $n\rightarrow\infty$ in the above inequality; then (\ref{45}) holds, and so (\ref{41}) holds.

By (\ref{41}) and the continuity of $\varphi$ in $\mathbb{R}^d$, we can show that for fixed $(t, s)\in\triangle$ and $x\in\mathbb{R}^d$,
\begin{equation}\nonumber
\begin{aligned}
\begin{split}
\varphi(t,s,\omega)[{\cal K}(u)](s,\omega)
& = \varphi(t,s,\omega) \lim \limits_{{n} \to {\infty}}{\varphi(s,-nT, \omega)x}\\
& = \lim \limits_{{n} \to {\infty}}\varphi(t,s,\omega){\varphi(s,-nT, \omega)x}\\
&=\lim \limits_{{n} \to {\infty}}{\varphi(t, -nT, \omega)x}\\
&=[{\cal K}(u)](t,\omega) \ \ \  \mathbb{P}-a.s.
\end{split}
\end{aligned}
\end{equation}
Furthermore, by (\ref{26}) and (\ref{44}) we have 
\begin{equation}\nonumber
\begin{aligned}
\begin{split}
[{\cal K}(u)](s,\theta_T\omega)
& = \displaystyle \int^{s}_{-\infty}{\Phi(s-r)u(r,\theta_T\omega)dr}+\displaystyle \int^{s}_{-\infty}{\Phi(s-r)\sigma(r) dW_r({\theta_T}\omega)}\\
&=\displaystyle \int^{s+T}_{-\infty}{\Phi(s+T-r)u(r-T,\theta_T\omega)dr}+\displaystyle \int^{s+T}_{-\infty}{\Phi(s+T-r)\sigma(r-T) dW_r}\\
&= \displaystyle \int^{s+T}_{-\infty}{\Phi(s+T-r)u(r,\omega)dr}+\displaystyle \int^{s+T}_{-\infty}{\Phi(s+T-r)\sigma(r) dW_r}\\
&= [{\cal K}(u)](s+T,\omega) \ \ \  \mathbb{P}-a.s.
\end{split}
\end{aligned}
\end{equation} 
\hfill\fbox\\

\begin{rem}
	Now we consider the nonlinear stochastic differential equation with the additive white noise of the form
	\begin{equation}\label{46}
	dX_t=[AX_t+h(X_t)]dt+\sigma dW_t, 
	\end{equation} 
	Here $h$ and $\sigma$ do not depend on time $t$, that is to say, the periodic $T$ in the Theorem 4.2 can be chosen as an arbitrary real number. The difference is that for this equation, we have a random dynamical system $\varphi: \mathbb{R}\times\Omega\times \mathbb{R}^d\rightarrow \mathbb{R}^d$, which satisfies the cocycle property $\varphi(t,\theta_s\omega)\varphi(s,\omega)=\varphi(t+s,\omega)$. Equation (\ref{46}) has a stationary solution. This result is given by J.F. Jiang and X. Lv ({\cite{JL}}).
\end{rem}

\section{Examples}
\quad In this section, we present several examples to illustrate the use of Theorem 4.2.
For simplicity, we only consider three-dimensional systems in the following. First, we give an example that $h$ is order-preserving in $\mathbb{R}^3$.
\begin{exa}
	Consider stochastic differential equation
	\begin{equation}\label{51}
	dx_i=[(Ax)_i+h_i(t,x_i)]dt+\sigma_i(t)dW_t^i, \ \ \ \ i=1,2,3,
	\end{equation}
	where
	$$A=\begin{bmatrix}
	-1&1&0\\
	1&-2&0\\
	0&1&-1
	\end{bmatrix}$$
	with three eigenvalues $\lambda_1=-1$, $\lambda_2=\frac{-3+\sqrt{5}}{2}$, $\lambda_3=\frac{-3-\sqrt{5}}{2}$,
	
	\begin{equation}\nonumber
	\begin{aligned}
	\begin{split}
	h_i(t,x_i)
	&=\frac{1}{10}sint+\frac{1}{6+g_i(x_i)}, \ \ i=1,2,3
	\end{split}
	\end{aligned}
	\end{equation}
	where 
	$g_i(x_i)=\frac{\pi}{2}-\arctan x_i$ is decreasing with respect to $x_i, i=1,2,3$ and $\sigma_i(t)=\cos t$. By direct calculation, we obtain 
	$$\Phi(t)=\begin{bmatrix}
	\frac{5+\sqrt{5}}{10}e^{\frac{-3+\sqrt{5}}{2}t}+\frac{5-\sqrt{5}}{10}e^{\frac{-3-\sqrt{5}}{2}t}&\frac{\sqrt{5}}{5}e^{\frac{-3+\sqrt{5}}{2}t}-\frac{\sqrt{5}}{10}e^{\frac{-3-\sqrt{5}}{2}t}&0\\
	\frac{\sqrt{5}}{5}e^{\frac{-3+\sqrt{5}}{2}t}-\frac{\sqrt{5}}{5}e^{\frac{-3-\sqrt{5}}{2}t}&\frac{5-\sqrt{5}}{10}e^{\frac{-3+\sqrt{5}}{2}t}+\frac{5+\sqrt{5}}{10}e^{\frac{-3-\sqrt{5}}{2}t}&0\\
	-e^{-t}+\frac{5+\sqrt{5}}{10}e^{\frac{-3+\sqrt{5}}{2}t}+\frac{5-\sqrt{5}}{10}e^{\frac{-3-\sqrt{5}}{2}t}&\frac{\sqrt{5}}{5}e^{\frac{-3+\sqrt{5}}{2}t}-\frac{\sqrt{5}}{5}e^{\frac{-3-\sqrt{5}}{2}t}&e^{-t}
	\end{bmatrix}$$
	It is not difficult to estimate that for any $t\geq 0$,
	$$||\Phi(t)||:= \max\{|\Phi_{ij}(t)|:i,j=1,2,3\}\leq e^{\lambda_2 t},$$
	which implies that (\ref{23}) holds. Moreover, it is easy to see 
	that $\max_{1\leq i\leq3}\text{Re} \lambda_i=\lambda_2\leq 0$, $L\leq \frac{1}{36}$. So
	$$-\frac{9L}{\lambda_2}\leq \frac{1}{2(3-\sqrt{5})}< 1.$$
	By Theorem 4.2, stochastic differential equation (\ref{51}) has a unique random periodic solution of periodic $2\pi$.
\end{exa}

Next, we give two examples that $h$ is anti-order-preserving in $\mathbb{R}^3$. 

\begin{exa}
	Consider stochastic differential equation
	\begin{equation}\label{52}
	dx_i=[a_ix_i+h_i(t,x_{i-1})]dt+\sigma_i(t)dW_t^i, \ \ \ \ i=1,2,3,
	\end{equation}
	where $a_1=-1, a_2=-2, a_3=-3$,  
	\begin{equation}\nonumber
	\begin{aligned}
	\begin{split}
	h_i(t,x_{i-1})
	&=\frac{1}{6+\cos t+{\rm th}\ x_{i-1}}:=\frac{1}{5+\cos t+g_i(x_{i-1})}, \ \ i=1,2,3
	\end{split}
	\end{aligned}
	\end{equation}
	where 
	$g_i(x_{i-1})=1+{\rm th}\ x_{i-1}$ is increasing with respect to $x_{i-1}, i=1,2,3$ and $\sigma_i(t)=\cos t$. It is easy to see 
	that $\lambda=-1$ and $L\leq \frac{1}{16}$, it follows that for any $t\geq 0$,
	$$||\Phi(t)||:= \max\{|\Phi_{ij}(t)|:i,j=1,2,3\}=e^{-t}$$
	and
	$$-\frac{9L}{\lambda}\leq \frac{9}{16}<1.$$
	By Theorem 4.2, stochastic differential equation (\ref{52}) has a unique random periodic solution of periodic $2\pi$.
\end{exa}

\begin{exa}
	Consider stochastic differential equation
	\begin{equation}\label{53}
	dx_i=[(Ax)_i+h_i(t,x_{i-1})]dt+\sigma_i(t)dW_t^i, \ \ \ \ i=1,2,3,
	\end{equation}
	where $x_0=x_3,x_4=x_1$ and  
	$$A=\begin{bmatrix}
	-1&\sqrt[3]{2}&0\\
	0&-2&\sqrt[3]{2}\\
	\sqrt[3]{2}&0&-4
	\end{bmatrix}$$
	with three eigenvalues $\lambda_1=-3$, $\lambda_2=-2+\sqrt{2}$, $\lambda_3=-2-\sqrt{2}$,
	
	\begin{equation}\nonumber
	\begin{aligned}
	\begin{split}
	h_i(t,x_{i-1})
	&:=\frac{1}{8}sint+\frac{1}{4+g_i(x_{i-1})}, \ \ i=1,2,3
	\end{split}
	\end{aligned}
	\end{equation}
	where 
	$g_i(x_{i-1})=\frac{\pi}{2}+\arctan x_{i-1}$ is increasing with respect to $x_{i-1}, i=1,2,3$ and $\sigma_i(t)=\sin t$. By direct calculation, we obtain 
	\begin{equation}\nonumber
	\begin{gathered}
	\begin{split}
	\Phi(t)=&\left[\begin{gathered}
	e^{-3t}+\frac{\sqrt{2}}{2}e^{(-2+\sqrt{2})t}-\frac{\sqrt{2}}{2}e^{-(2+\sqrt{2})t}\\	-\sqrt[3]{4}e^{-3t}+\frac{\sqrt[3]{4}-\sqrt[6]{2}}{2}e^{(-2+\sqrt{2})t}+\frac{\sqrt[3]{4}+\sqrt[6]{2}}{2}e^{-(2+\sqrt{2})t}\\\sqrt[3]{2}e^{-3t}+\frac{\sqrt[6]{2^5}-\sqrt[3]{2}}{2}e^{(-2+\sqrt{2})t}-\frac{\sqrt[6]{2^5}+\sqrt[3]{2}}{2}e^{-(2+\sqrt{2})t}
	\end{gathered}\right.&\\
	&\begin{gathered}
	-\sqrt[3]{2}e^{-3t}+\frac{\sqrt[3]{2}}{2}e^{(-2+\sqrt{2})t}+\frac{\sqrt[3]{2}}{2}e^{-(2+\sqrt{2})t}\\2e^{-3t}+\frac{\sqrt{2}-1}{2}e^{(-2+\sqrt{2})t}-\frac{\sqrt{2}+1}{2}e^{-(2+\sqrt{2})t}\\-\sqrt[3]{4}e^{-3t}+\frac{\sqrt[3]{4}-\sqrt[6]{2}}{2}e^{(-2+\sqrt{2})t}+\frac{\sqrt[3]{4}+\sqrt[6]{2}}{2}e^{-(2+\sqrt{2})t}
	\end{gathered}&\\
	&\left.\begin{gathered}
	-\sqrt[3]{4}e^{-3t}+\frac{\sqrt[3]{4}-\sqrt[6]{2}}{2}e^{(-2+\sqrt{2})t}+\frac{\sqrt[3]{4}+\sqrt[6]{2}}{2}e^{-(2+\sqrt{2})t}\\2\sqrt[3]{2}e^{-3t}+(\frac{3\sqrt[6]{2^5}}{4}-\sqrt[3]{2})e^{(-2+\sqrt{2})t}-(\frac{3\sqrt[6]{2^5}}{4}+\sqrt[3]{2})e^{-(2+\sqrt{2})t}\\-2e^{-3t}+(\frac{3}{2}-\sqrt{2})e^{(-2+\sqrt{2})t}+(\frac{3}{2}+\sqrt{2})e^{-(2+\sqrt{2})t}
	\end{gathered}\right]&,
	\end{split}
	\end{gathered}
	\end{equation}
	It is not difficult to prove that for any $t\geq 0$,
	$$||\Phi(t)||:= \max\{|\Phi_{ij}(t)|:i,j=1,2,3\}\leq e^{\lambda_2 t},$$
	which can be found in ({\cite{JL}}, Example 5.3) for a detailed proof. Moreover, it is easy to see 
	that $\max_{1\leq i\leq3}\text{Re} \lambda_i=\lambda_2\leq 0$, $L\leq \frac{1}{16}$, and so
	$$-\frac{9L}{\lambda_2}\leq \frac{9}{16(2-\sqrt{2})}< 1.$$
	By Theorem 4.2, stochastic differential equation (\ref{53}) has a unique random periodic solution of periodic $2\pi$.
\end{exa}

\vskip 1cm

{ \noindent {\bf\large Acknowledgements}\ \
	The authors would like to express his deep gratitude to Professor Jifa Jiang, who gave numerous suggestions which significantly improve the readability of the paper. This work was supported in part by National Natural Science Foundation of China (Nos. 11431014, 11671382, 11931004, 12031020, 12090014), CAS Key Project of Frontier Sciences (No. QYZDJ-SSW-JSC003), the Key Laboratory of Random Complex Structures and Data Sciences, CAS(No. 2008DP173182) and National Center for Mathematics and Interdisciplinary Sciences CAS.


\begin{thebibliography}{1234} 
		\bibitem{A} L. Arnold, Random Dynamical Systems, Springer Monogr. Math., Springer-Verlag, Berlin, 1998, doi:10.1007/978-3-662-11478-7.	
		
		\bibitem{C} I. Chueshov, Monotone Random Systems Theory and Applications, Lecture Notes in Math. 1799, Springer-Verlag, Berlin, 2002, doi:10.1007/b83277. 	
		
		\bibitem{FLZ} C.R. Feng, Y. Liu and H.Z. Zhao, Numerical approximation of random periodic solutions of stochastic differential equations, Z. Angew. Math. Phys. 68(2017) 119.
		
		\bibitem{FLZZ} C.R. Feng, Y. Liu and H.Z. Zhao, Numerical analysis of the weak schemes of random periodic solutions of stochastic differential equations, submitted for publication.
		
		\bibitem{FWZ} C.R. Feng, Y. Wu and H.Z. Zhao, Anticipating random periodic solutions-I. SDEs with multiplicative linear noise, J. Funct. Anal. {271}, 365-417 (2016).
		
		\bibitem{FZ} C.R. Feng and H.Z. Zhao, Random periodic solutions of SPDEs via integral equations and Wiener-Sobolev compact embedding, J. Funct. Anal. {262}, 4377-4422 (2012).
		
		\bibitem{FFZ} C.R. Feng and H.Z. Zhao, Random periodic processes, periodic measures and ergodicity, J. Differ. Equ. {269}, 7382-7414 (2020).
		
		\bibitem{FZZ} C.R. Feng, H.Z. Zhao and B. Zhou, Pathwise random periodic solutions of stochastic differential equations, J. Differ. Equ. {251}, 119-149 (2011).
		
		\bibitem{JL} J.F. Jiang and X. Lv, A small-gain theorem for nonlinear stochastic systems with inputs and outputs I: Additive white noise, SIAM J. Control Optim., 54(2016), pp. 2383-2402. 
		
		\bibitem{KS} I. Karatzas and S.E. Shreve, Brownian Motion and Stochastic Calculus, Grad. Texts in Math. 113, Springer, New York, 1988, doi:10.1007/978-1-4612-0949-2. 
		
    	\bibitem{K} H. Kunita, 1984, Stochastic Differential Equations and Stochastic Flows of Diffeomorphisms, Ecole d'\'et\'e de Probabilit\'es de Saint-Flour 12, 1982 Lect.Notes Math.1097, 143-303.  	   	
				
		\bibitem{KH} H. Kunita, Stochastic Flows and Stochastic Differential Equations, Cambridge University Press, Cambridge(1990).     	   	
		
		\bibitem{MS} M. Marcondes de Freitas and E.D. Sontag, A small-gain theorem for random dynamical systems with inputs and outputs, SIAM J. Control Optim., 53(2015), pp. 2657-2695, doi:10.1137/140991340.    
		
		\bibitem{M} X.R. Mao, Stochastic Differential Equations and Applications, Horwood, Chichester, UK, 1997.      	
		
		\bibitem{O} B.\O. ksendal, Stochastic Differential Equations: An Introduction with Applications, 5th ed., Springer-Verlag, Berlin, 1998, doi:10.1007/978-3-642-14394-6. 
		
		\bibitem{PM} J. Palis and W. Melo, Geometric Theory of Dynamical Systems, An Introduction, Springer-Verlag, New York, 1982, doi:10.1007/978-1-4612-5703-5. 
		
		\bibitem{P} H. Poincar\'e, memoire sur les courbes definier par une equation differentitate, J.Math.Pures Appl.3(1881) 375-442, J.Math.Pures Appl.3(1882) 251-296, J.Math.Pures Appl.4(1885) 167-244, J.Math.Pures Appl.4(1886) 151-217.    	       	
		
		\bibitem{Y} K. Yoshida, Functional Analysis, 6th ed., Springer, New York, 1980, doi:10.1007/978-3-642-61859-8.  	
		
		\bibitem{ZZ} H.Z. Zhao and Z-H. Zheng, Random periodic solutions of random dynamical systems, J. Differ. Equ. {246}, 2020-2038 (2009). 
		
		
	\end{thebibliography}
\end{document}